\documentclass[11pt]{article}

\usepackage[draft]{graphics}
\usepackage{graphicx}          	 
\usepackage{bm}                 
\usepackage{amsmath}             
\usepackage{amssymb}
\usepackage{amsfonts}             
\usepackage{verbatim}           
\usepackage{amsthm}             

\usepackage{mathtools}
\usepackage{relsize}

\usepackage{makeidx}

\theoremstyle{plain}             

\theoremstyle{definition}

\makeatletter
\@addtoreset{equation}{chapter}

\makeatother

\def\Fgauss#1#2#3#4{
{}_2F_1\left(
\begin{array}{c}
\begin{array}{cc} \hskip-10pt#1,{\ } #2 \end{array}\\
\begin{array}{c} \hskip-10pt#3 \end{array}
\end{array}
\hskip-8pt;\,#4
\right)}

\def\protectbold#1{\protect{\boldmath{$#1$}}}

\def\eqref#1{(\ref{#1})}
\def\dsp{\displaystyle}
\def\Frac#1#2{\frac
{
 {\raise.6ex
 \hbox{$\displaystyle#1$}}
}
{
 {\lower.6ex
 \hbox{$\displaystyle#2$}}
 }
}

\numberwithin{equation}{section}

\def\bigOxe{\sqcup \kern-2.3mm \sqcap}

\makeatother

\def\dsp{\displaystyle}
\def\Frac#1#2{\frac
{
 {\raise.6ex
 \hbox{$\displaystyle#1$}}
}
{
 {\lower.6ex
 \hbox{$\displaystyle#2$}}
 }
}

\input epsf

\def\CHFs#1#2#3{
{}_1F_1\left({a};{c};{z}\right)
}

\def\FG#1#2#3#4{
{}_2F_1\left(
\begin{array}{c}
\begin{array}{c}\hskip-10pt#1,#2\end{array}\\
\begin{array}{c}\hskip-10pt #3\end{array}
\end{array}
\hskip-8pt;\,#4
\right)}

\def\intp{\int_0^\infty}

\def\bigO{{\cal O}}

\def\tfrac#1#2{{{\lower.6ex
\hbox{$\scriptstyle#1$}}\over 
{\raise.7ex
\hbox{$\scriptstyle#2$}}}}

\def\intp{\int_0^\infty}

\def\phase{{\rm ph}}

\def\tfrac#1#2{{{\lower.6ex
\hbox{$\scriptstyle#1$}}\over 
{\raise.7ex
\hbox{$\scriptstyle#2$}}}}

\def\sn{\sum_{n=0}^\infty\,}
\def\sk{\sum_{k=0}^\infty\,}

\def\insil#1{}

\begin{document}
 \title{Computation of a numerically satisfactory pair of solutions of the differential
equation for conical functions of non-negative integer orders}

\author{
T.M. Dunster\\
Department of Mathematics and Statistics\\
San Diego State University. 5500 Campanile Drive San Diego, CA, USA.
\and
A. Gil\\
Departamento de Matem\'atica Aplicada y CC. de la Computaci\'on.\\
ETSI Caminos. Universidad de Cantabria. 39005-Santander, Spain.\\
\and
J. Segura\\
        Departamento de Matem\'aticas, Estad\'{\i}stica y 
        Computaci\'on,\\
        Univ. de Cantabria, 39005 Santander, Spain.\\   
    \and
    N.M. Temme\\
    IAA, 1391 VD 18, Abcoude, The Netherlands\footnote{Former address: Centrum Wiskunde \& Informatica (CWI), 
        Science Park 123, 1098 XG Amsterdam,  The Netherlands}\\
}

\date{\ }

\maketitle
\begin{abstract}
We consider the problem of computing satisfactory pairs of solutions of the differential equation
for Legendre functions of non-negative integer order $\mu$ and degree $-\frac12+i\tau$, where $\tau$ is a non-negative
 real parameter. 
Solutions of this equation are the conical functions ${\rm{P}}^{\mu}_{-\frac12+i\tau}(x)$ and 
${Q}^{\mu}_{-\frac12+i\tau}(x)$, $x>-1$.
An algorithm for computing a numerically satisfactory pair of solutions 
 is already available when $-1<x<1$ (see \cite{gil:2009:con}, \cite{gil:2012:cpc}).
In this paper, we present a stable computational scheme for
 a real valued numerically satisfactory companion of the function ${\rm{P}}^{\mu}_{-\frac12+i\tau}(x)$
for $x>1$, the function $\Re\left\{e^{-i\pi \mu} {{Q}}^{\mu}_{-\frac{1}{2}+i\tau}(x)    \right\}$. 
 The proposed algorithm allows the computation of the function 
on a large parameter domain without requiring the use of extended precision arithmetic.
\end{abstract}

\section{Introduction}

Conical or Mehler functions appear in a large number of applications in applied physics, 
particle physics (related to the amplitude for Yukawa potential 
scattering)
or cosmology, among others. These functions also appear when solving the Laplace's
problem in spherical coordinates for two intersecting cones \cite{thebault:2006:geo}.  

  The conical functions are solutions of the second order differential equation

\begin{equation}\label{ODE}
(1-x^2) \Frac{d^2 w}{dx^2}-2x \Frac{dw}{dx}-\left(\tau^2+\Frac{1}{4}+\Frac{\mu^2}{1-x^2}  \right)w=0.
\end{equation}
   
 In applications, one encounters values of the conical functions where
the order $\mu$ is integer positive (or zero), so we will consider $\mu=m \in {\mathbb Z}^+$.

\section{Numerically satisfactory pairs of solutions} 

In the interval $-1 <x <1$, a numerically satisfactory pair of solutions (real valued) of the equation
\eqref{ODE}
is ${\rm P}_{-\frac12+i\tau}^{-\mu}(x)$ and ${\rm P}_{-\frac12+i\tau}^{-\mu}(-x)$, $\mu \in \Re$, $\tau \in \Re^+$  \cite{Duns:2010:Con}.  The Wronskian
for this pair of solutions is:

\begin{equation}
\label{W1}
W\left\{{\rm{P}}^{-\mu}_{-\frac{1}{2}+i\tau}\left(x\right),{\rm{P}}^{-\mu}_{-\frac{1}{2}+i\tau}\left(-x\right)\right\}
=\Frac{2}{\left|\Gamma\left(\mu+\tfrac{1}{2}+i\tau\right)\right|^{2}(1-x^2)}.
\end{equation}

When $\mu=m \in {\mathbb Z}$, it should be noted that the 
functions ${\rm P}_{-\frac12+i\tau}^{m}(x)$ and ${\rm P}_{-\frac12+i\tau}^{-m}(x)$ are related
  through a simple relation:

\begin{equation}
\label{changem}
{\rm{P}}^m_{-\frac12+i\tau}(x)=\cosh(\pi\tau)\Frac{\left|\Gamma (m+\tfrac12+i\tau)\right|^2}{\pi}
{\rm{P}}^{-m}_{-\frac12+i\tau}(x).
\end{equation}

The definition of the function ${\rm P}^m_{-\frac12+i\tau}(x)$ in terms of the Gauss 
hypergeometric function $_2F_1$ is given in Eq. \eqref{changem}.
 In \cite{gil:2009:con} and \cite{gil:2012:cpc} we have described
an algorithm for computing the conical functions ${\rm P}_{-\frac12+i\tau}^{m}(x)$ for  
$x>-1$; therefore the problem of computing a numerically satisfactory pair of solutions
of the differential equation for conical functions in the interval $-1<x<1$
can be considered as already solved. The algorithm is based on the use of different methods of computation, depending on the range
of the parameters:
quadrature methods, recurrence relations and  uniform 
asymptotic expansions in terms of elementary functions or in terms of 
modified Bessel function $K_{ia}(x)$ and its derivative $K'_{ia}(x)$.  

\begin{figure}
\caption{Numerically satisfactory pair of solutions of \eqref{ODE} for $x>1$.
The functions ${\rm P}^{10}_{-\frac{1}{2}+i5}(x)$ and  $\widetilde{{Q}}^{10}_{-\frac{1}{2}+i5}(x)$
are plotted as an example. 
\label{Fig1}}
\begin{center}
\epsfxsize=13cm \epsfbox{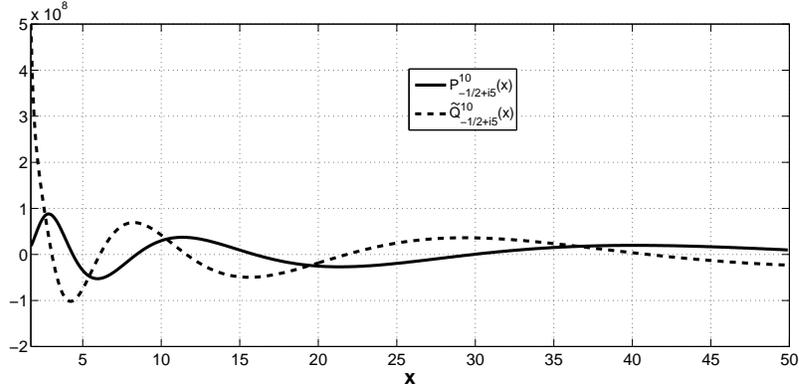}
\end{center}
\end{figure}

 When $x>1$, a real-valued satisfactory pair of solutions is
${\rm P}^{m}_{-\frac{1}{2}+i\tau}(x)$ and $\widetilde{{Q}}^{m}_{-\frac{1}{2}+i\tau}(x) 
\equiv \Re\left\{e^{-i\pi m} {{Q}}^{m}_{-\frac{1}{2}+i\tau}(x)    \right\}$
(the function ${{Q}}^{m}_{-\frac{1}{2}+i\tau}(x)$ is complex-valued). The solution  $\widetilde{{Q}}^{m}_{-\frac{1}{2}+i\tau}(x)$ was introduced
in \cite{dunster:1991:cfo}. As an example, Fig.~\ref{Fig1} shows the graphs of
${\rm P}^{10}_{-\frac{1}{2}+i5}(x)$ and  $\widetilde{{Q}}^{10}_{-\frac{1}{2}+i5}(x)$.

 As in the previous case, an algorithm for computing the function ${\rm P}^{m}_{-\frac{1}{2}+i\tau}(x)$ for $x>1$ is already
given in \cite{gil:2009:con} and \cite{gil:2012:cpc}. However, no computational schemes can be found in the literature for computing
the numerically satisfactory companion solution of ${\rm P}^{m}_{-\frac{1}{2}+i\tau}(x)$, the function  $\widetilde{{Q}}^{m}_{-\frac{1}{2}+i\tau}(x)$.
Our purpose is to give an algorithm for computing this function. 

\section{Recurrence relations}\label{RR}

Three-term recurrence relations 
\begin{equation}
y_{n+1}+b_n y_n +a_n y_{n-1}=0 ,
\end{equation}
are useful methods of computation when
two starting values are available for starting the recursive process. 
Usually, 
the direction of application of the recursion can not be chosen arbitrarily,
and the conditioning of the computation of a given solution
fixes the direction. 

The conical function ${\rm P}_{-1/2+i\tau}^{m}(x)$ satisfies a three-term recurrence relation,
\begin{equation}
\label{TTRR2}
{\rm P}_{-\frac12+i\tau}^{m +1}(x)-\Frac{2m x}{\sqrt{x^2-1}}
{\rm P}_{-\frac12+i\tau}^{m}(x) +
\left((m-\tfrac12)^2+\tau^2\right){\rm P}_{-\frac12+i\tau}^{m-1}(x)  =0,
\end{equation}
for $x>1$, where we have adopted the following definition of the function ${\rm P}^m_{-\frac12+i\tau}(x)$ in terms of the Gauss 
hypergeometric function $_2F_1$:

\begin{equation}
\label{changem}
\begin{array}{lcl}
{\rm P}^m_{-\frac12+i\tau}(x)&=&\cosh(\pi\tau)\Frac{|\Gamma (m+1/2+i\tau)|^2}{\pi \Gamma(1+m)}
\left|\Frac{1-x}{1+x}\right|^{m/2} \times \\ 
&&\Fgauss{\tfrac12 -i \tau}{\tfrac12 +i \tau }{1+m}{\tfrac12-\tfrac12x}.
\end{array}
\end{equation}

This is the definition used in \cite{gil:2012:cpc} for $x>-1$.
Note the difference of the sign in the second term of the recurrence relation \eqref{TTRR2} in comparison with
Eq. 14.10.6 in \cite{Duns:2010:Con}, which reflects the fact that the definition of our function ${\rm P}^m_{-\frac12+i\tau}(x)$
and the function $P_{-1/2+\i\tau }^{m}\left ( x \right )$ (\cite{dunster:1991:cfo}, \cite{Duns:2010:Con})  differ
by a factor $(-1)^m$:

  \begin{equation}
   \textup{P}_{-1/2+\i\tau  }^{m}\left ( x \right )=\left ( -1 \right )^{m}P_{-1/2+\i\tau }^{m}\left ( x \right )\quad\ (x>1).
  \end{equation}

The conical functions ${\rm P}_{-\frac12+i\tau}^{m}(x)$ are monotonic in the interval $(1,x_c)$ and oscillating in $(x_c,+\infty)$, where
$x_c=\sqrt{1+\beta^2}/\beta$ and $\beta=\tau/m$. As an example, for the parameters of the functions plotted in Fig.~\ref{Fig1}, the
transition point is found at $x_c \approx 2.24$.
The stability analysis based on Perron's theorem discussed in \cite{gil:2009:con}, revealed that
the function  ${\rm P}_{-\frac12+i\tau}^{m}(x)$ was the minimal solution of the recurrence relation \eqref{TTRR2}. Hence,
for this function backward recursion will be 
generally stable for $x>1$. However, and similar for other special functions, recurrence relations in 
the oscillatory regime of the conical functions ($x>x_c>1$)
are not ill conditioned in either backward and forward directions; thus, both recursions are possible.

It is very simple to show that the function $\widetilde{{Q}}^{m}_{-\frac{1}{2}+i\tau}(x)$ also satisfies 
the three-term recurrence relation
given in \eqref{TTRR2}, and this function is a dominant solution of \eqref{TTRR2} (the minimal solution, if it exists,
is unique).
In this case, the stable direction of application of the recurrence relation is with increasing $m$, although the same comment
as for the function ${\rm P}_{-\frac12+i\tau}^{m}(x)$ applies for  $\widetilde{{Q}}^{m}_{-\frac{1}{2}+i\tau}(x)$ in the
oscillatory regime ($x>x_c>1$).

 We next consider the problem of computing two starting values ($\widetilde{{Q}}^{0}_{-\frac{1}{2}+i\tau}(x)$
and  $\widetilde{{Q}}^{1}_{-\frac{1}{2}+i\tau}(x)$) of the recurrence  relation \eqref{TTRR2} for the
function  $\widetilde{{Q}}^{m}_{-\frac{1}{2}+i\tau}(x)$.

\section{Power series expansions by using the hypergeometric functions}\label{sec:pow1}
From the many representations of the conical functions in terms of the Gauss hypergeometric functions we choose one representation that 
can be used for small values of $x-1$ and one for large values of $x$ (and that one can also be used for large values of $\tau$).

\subsection{Expansions valid near the point \protectbold{x=1} and moderate $\tau$}\label{sec:pow2}
We use the representation (see \cite[Eqn.~(32), page 131]{Bateman:1953:HTF})
\begin{equation}\label{eq:QKumF01}
e^{-\pi i\mu}{Q}_{\nu}^\mu(x)=
A_3\,\FG{-\nu}{1+ \nu}{1+\mu}{z}
+A_4\,\FG{-\nu}{1+ \nu}{1-\mu}{z},
\end{equation}
where
\begin{equation}\label{eq:QKumF02}
\begin{array}{@{}r@{\;}c@{\;}l@{}@{}r@{\;}c@{\;}l@{}}
A_3&=&\dsp{\tfrac12\Gamma(-\mu)w^{\mu}\frac{\Gamma(1+\nu+\mu)}{\Gamma(1+\nu-\mu)}},
\quad\quad &A_4=&\dsp{\tfrac12\Gamma(\mu)w^{-\mu},}&\\[8pt]
z&=&\tfrac12(1-x), &w=\hskip-8pt&\dsp{\sqrt{\frac{x-1}{x+1}}.}&
\end{array}
\end{equation}
We use power series in $z$, and because we want these expansion for $\mu=0$ and $\mu=1$, we need to find the limit for the representation in \eqref{eq:QKumF01}. 

First we write
\begin{equation}\label{eq:QKumF03}
2\Gamma(1+\nu-\mu)e^{-\pi i\mu}{Q}_{\nu}^\mu(x)=\frac{\pi \mu}{\sin(\pi\mu)}\sk \frac{(-\nu)_k(1+\nu)_k}{k!}z^kB_k(\mu),
\end{equation}
where
\begin{equation}\label{eq:QKumF04}
B_k(\mu)=\frac{1}{\mu}\left(w^{-\mu}\frac{\Gamma(1+\nu-\mu)}{\Gamma(1-\mu+k)}-w^{\mu}\frac{\Gamma(1+\nu+\mu)}{\Gamma(1+\mu+k)}\right).
\end{equation}
The limit for $\mu\to0$ gives
\begin{equation}\label{eq:QKumF05}
B_k(0)=\frac{2\Gamma(\nu+1)}{k!}\bigl(\psi(k+1)-\psi(1+\nu)-\ln w\bigr),
\end{equation}
where $\psi(\alpha)=\Gamma^\prime(\alpha)/\Gamma(\alpha)$. We obtain the expansion
\begin{equation}\label{eq:QKumF06}
{Q}_{\nu}^{\,0}(x)=\sk \frac{(-\nu)_k(1+\nu)_k}{k!\,k!}z^k\bigl(\psi(k+1)-\psi(1+\nu)-\ln w\bigr),
\end{equation}
where $w$ and $z$ are given in \eqref{eq:QKumF02}.

Observe that for $\nu=-\frac12+i\tau$ we have
\begin{equation}\label{eq:QKumF07}
(-\nu)_k(1+\nu)_k=\left(\tfrac12-i\tau\right)_k\left(\tfrac12+i\tau\right)_k,
\end{equation}
which is real when $\tau$ is real, and the computation easily follows from recursion. The quantity $\psi(\nu+1)$ is the only complex term to consider in more detail. For $\psi(k+1)$ we can use the recursion $\psi(k+1)=\psi(k)+1/k$, with $\psi(1)=-\gamma$ ($\gamma$ is Euler's constant).

For the computation of $\psi(\alpha)$ (with $\alpha$ complex) we use the asymptotic expansion
\begin{equation}\label{eq:QKumF08}
\begin{array}{@{}r@{\;}c@{\;}l@{}}
\psi(\alpha)&\sim&\dsp{\ln \alpha-\frac{1}{2\alpha}-\sum_{n=1}^\infty \frac{B_{2n}}{2n\,\alpha^{2n}}}\\[8pt]
&=&\dsp{\ln \alpha}-\frac{1}{2\alpha}-\frac{1}{12\alpha^2}+\frac{1}{120\alpha^4}-\frac{1}{252\alpha^6}+...,
\end{array}
\end{equation}
valid for $\alpha\to\infty$ in $\vert\phase\,\alpha\vert<\pi$. We use this expansion if $\vert \alpha\vert \ge 12$ with $8$ terms 
of the series (or less), and we use the backward recurrence relation $\psi(\alpha)=\psi(\alpha+1)-1/\alpha$ 
for smaller values of $\vert \alpha\vert$. We only need an algorithm for $\Re \alpha\ge	\frac12$.

For ${Q}_{\nu}^{\,1}(x)$ we can also use a limiting procedure, but it is more convenient to use the derivative of ${Q}_{\nu}^{\,0}(x)$, because
\begin{equation}\label{eq:QKumF09}
{Q}_{\nu}^{\,1}(x)=\sqrt{x^2-1}\frac{d}{dx}{Q}_{\nu}^{\,0}(x).
\end{equation}
Expansions for integer values of $\mu$ are also given in \cite[Eqn.~(32), page 131]{Bateman:1953:HTF}).

\subsection{Expansions for moderate or large values of  \protectbold{x} and \protectbold{\tau}}\label{sec:pow3}
For moderate or large values of $x$ and $\tau$ we use  the representation  (see \cite[\S3.6.1, page 149]{Bateman:1953:HTF})
\begin{equation}\label{eq:QKum01}
e^{-\pi i\mu}{Q}_{-\frac12+i\tau}^\mu(x)=
A\,\FG{\frac12+\mu}{\frac12-\mu}{1+i\tau}{-z},
\end{equation}
where
\begin{equation}\label{eq:QKum02}
\begin{array}{@{}r@{\;}c@{\;}l@{}}
A&=&\dsp{\sqrt{\pi/2}\left(x^2-1\right)^{-\frac14}\left(x+\sqrt{x^2-1}\right)^{-i\tau}
\frac{\Gamma\left(\mu+i\tau+\frac12\right)}{\Gamma\left(1+i\tau\right)},}\\[8pt]
z&=&\dsp{\frac{1}{2\sqrt{x^2-1}\left(x+\sqrt{x^2-1}\right)}.}
\end{array}
\end{equation}
This gives the expansion

\begin{equation}\label{eq:hyp1}
\begin{array}{@{}r@{\;}c@{\;}l@{}}
e^{-\pi i\mu}{Q}_{-\frac12+i\tau}^\mu(x) &=&\dsp{\sqrt{\pi/2}\left(x^2-1\right)^{-\frac14}\left(x+\sqrt{x^2-1}\right)^{-i\tau}}
\ \times \\
&&\dsp{\frac{\Gamma\left(\tfrac12+\mu+i\tau\right)}{\Gamma\left(1+i\tau\right)}\sum_{k=0}^{\infty} \Frac{\left(\tfrac12+\mu \right)_k\left(\tfrac12-\mu \right)_k }
{\left(1+i\tau\right)_k}\Frac{(-z)^k}{k!}}.
\end{array}
\end{equation}

When computing the expansion for $\tau$ large, it is convenient to use the
following asymptotic expansion for the ratio of two gamma functions \cite{Temme:1996:SFA}:

\begin{equation}
\label{ratgam}
\Frac{\Gamma(z+a)}{\Gamma(z+b)} \sim z^{a-b} \displaystyle\sum_{n=0}^{\infty} c_n 
\Frac{\Gamma (b-a+n)}{\Gamma (b-a)}\Frac{1}{z^{n}},\,\,\mbox{as } z \rightarrow \infty,
\end{equation}
in the sector $|\arg z| <\pi$, with $a$ and $b$ fixed.

The coefficients $c_n$ appearing in \eqref{ratgam} are given in terms of generalized Bernoulli polynomials by

\begin{equation}
c_n=(-1)^n\Frac{B_n^{(a-b+1)}(a)}{n!}.
\end{equation}
In the present case, we have
\begin{equation}
\label{ratgam2}
\Frac{\Gamma(\tfrac12+\mu+i\tau)}{\Gamma(1+i\tau)} \sim \tau^{\mu-\frac12}e^{\frac12\left(\mu-\tfrac12\right)\pi i} \displaystyle\sum_{n=0}^{\infty} (-i)^n  \left(\tfrac12-\mu\right)_n\frac{c_n}{\tau^n}, 
\end{equation}
and the first coefficients are
\begin{equation}
\label{ratgam3}
c_0=1, \quad c_1=-\tfrac14(2\mu+1),\quad c_2=\tfrac1{96}(2\mu+1)(6\mu+1).
\end{equation}

The representation in \eqref{eq:hyp1} becomes
\begin{equation}\label{eq:qhyp1}
\begin{array}{@{}r@{\;}c@{\;}l@{}}
e^{-\pi i\mu}{Q}_{-\frac12+i\tau}^\mu(x)&=&
 \dsp{\sqrt{\pi/2} \,\tau^{\mu-\frac12}\left(x^2-1\right)^{-\frac14}G(\mu,\tau)\,e^{-i\phi}}
\ \times \\
&&\dsp{\sum_{k=0}^{\infty} \left(\tfrac12+\mu \right)_k\left(\tfrac12-\mu \right)_k 
\frac{u_k(\tau)+iv_k(\tau)}{w_k(\tau)}\Frac{(-z)^k}{k!},                                           }
\end{array}
\end{equation}
where $\phi=\tau\log\left(x+\sqrt{x^2-1}\right)-\frac12\left(\mu-\frac12\right)\pi$, 

\begin{equation}
\label{gmu}
G(\mu,\tau)=\displaystyle\sum_{n=0}^{\infty} (-i)^n  \left(\tfrac12-\mu\right)_n\frac{c_n}{\tau^n}, 
\end{equation}
and
\begin{equation}
\label{eq:qhyp2}
\frac{u_k(\tau)+iv_k(\tau)}{w_k(\tau)}=\frac{1}{\left(1+i\tau\right)_k}, \quad k=0,1,2,\ldots.
\end{equation}
We can obtain these quantities from the recurrence relations
\begin{equation}
\label{eq:qhyp3}
\begin{array}{@{}r@{\;}c@{\;}l@{}}
u_{k+1}(\tau)&=&(k+1)u_k(\tau)+\tau v_k(\tau),\\[8pt]
v_{k+1}(\tau)&=&(k+1)v_k(\tau)-\tau u_k(\tau),\\[8pt]
w_{k+1}(\tau)&=&\left((k+1)^2+\tau^2\right)w_k(\tau),
\end{array}
\end{equation}
with $u_0(\tau)=1$, $v_0(\tau)=0$, $w_0(\tau)=1$.

We can write \eqref{eq:hyp1} with a single complex representation, by writing
\begin{equation}
\label{eq:qhyp4}
G(\mu,\tau)=H(\mu,\tau)e^{i\rho(\mu,\tau)},\quad u_k(\tau)+iv_k(\tau)=r_k(\tau)e^{i\sigma_k(\tau)},
\end{equation}
which gives
\begin{equation}
\label{eq:qhyp5}
\begin{array}{@{}r@{\;}c@{\;}l@{}}
e^{-\pi i\mu}{Q}_{-\frac12+i\tau}^\mu(x)&=&
 \dsp{\sqrt{\pi/2}\, H(\mu,\tau)\,\tau^{\mu-\frac12}\left(x^2-1\right)^{-\frac14}\,
\ \times} \\
&&\dsp{\sum_{k=0}^{\infty} \left(\tfrac12+\mu \right)_k\left(\tfrac12-\mu \right)_k 
\Frac{r_k(\tau)}{w_k(\tau)}\Frac{(-z)^k}{k!}e^{-i\psi_k},                                           }
\end{array}
\end{equation}
where
\begin{equation}
\label{eq:qhyp6}
\psi_k =\tau\log\left(x+\sqrt{x^2-1}\right)-\tfrac12\left(\mu-\tfrac12\right)\pi-\rho(\mu,\tau)-\sigma_k(\tau).
\end{equation}

\section{Expansions in terms of the Kummer \protectbold{U-}functions}\label{sec:QasyKum1}

First we explain how the method works for the Gauss hypergeometric function.
We take the integral representation 
\begin{equation}\label{eq:QKum03}
\FG{a}{b}{c}{z}=\frac{\Gamma(c)}{\mathop{\Gamma\/}\nolimits\!\left(b\right)\mathop{\Gamma\/}\nolimits\!\left(c-b\right)}\int _{0}^{1}u^{{b-1}}(1-u)^{{c-b-1}}(1-z u)^{-a}\,du, 
\end{equation}
and using the transformation $u=1-e^{-t}$, we write it in the form
\begin{equation}\label{eq:QKum04}
(1+z)^a\frac{\Gamma(c+\omega-b)}{\Gamma(c+\omega)}
 \FG{a}{b}{c+\omega}{-z}=F_{a,b}(\alpha,\omega),
\end{equation}
where
\begin{equation}\label{eq:QKum05}
F_{a,b}(\alpha,\omega)=\frac{1}{\Gamma(b)}\intp t^{b-1}f(t)e^{-\omega t}\,\frac{dt}{(t+\alpha)^a}, 
\end{equation}
with
\begin{equation}\label{eq:QKum06}
f(t)=\left(\frac{e^t-1}{t}\right)^{b-1}e^{(1+a-c)t}\left(\frac{e^t-e^{-\alpha}}{t+\alpha}\right)^{-a},\quad \alpha=\ln\frac{z+1}{z}.
\end{equation}

We assume that $\omega$ is large and that  $z$ may be large as well. In that case $\alpha$ will be small.
The easiest approach is to expand
$\dsp{f(t)=\sn f_nt^n}$, which gives
\begin{equation}\label{eq:QKum07}
F_{a,b}(\alpha,\omega)\sim \sn f_n \Phi_n,
\end{equation}
where
\begin{equation}\label{eq:QKum08}
\Phi_n=\frac1{\Gamma(b)}\intp
t^{n+b-1}(t+\alpha)^{-a}e^{-\omega t}\,dt .
\end{equation}

The functions $\Phi_n$ can be expressed in terms of the confluent hypergeometric function $U(a,c,z)$. We have
\begin{equation}\label{eq:QKum09}
\begin{array}{@{}r@{\;}c@{\;}l@{}}
\Phi_n &=&
\dsp{(b)_n \alpha^{n+b-a}U(n+b,n+b+1-a,\alpha \omega)}\\[8pt]
&=&
\dsp{(b)_n\omega^{a-n-b}U(a,a+1-n-b,\alpha \omega).}
\end{array}
\end{equation}

The functions $\Phi_n$ can be obtained by using recurrence relations for the $U-$function\footnote{http://dlmf.nist.gov/13.3.E14}, or from integrating by parts in  \eqref{eq:QKum08}.
We have
\begin{equation}\label{eq:QKum10}
\omega\Phi_{n+1}=(n+b-a-\alpha\omega)\Phi_n+\alpha(b+n-1)\Phi_{n-1}.
\end{equation}

For numerical aspects of such recursions, we refer to \cite{Temme:1983:NCC} and \cite[\S4.5.1]{Gil:2007:NSF}.

The expansion in \eqref{eq:QKum07} has an asymptotic character for large values of $\omega$, uniformly with respect to $\alpha$. The second line of \eqref{eq:QKum09} gives the integral representation
\begin{equation}\label{eq:QKum11}
\Phi_n=\frac{(b)_n}{\omega^{n+b-a}\Gamma(a)} \intp t^{a-1}(1+t)^{-n-b}e^{-\alpha \omega t}\,dt.
\end{equation}
which shows that 
\begin{equation}\label{eq:QKum12}
\Phi_n=\bigO\left(\omega^{-n-b+a}\right),\quad \omega\to\infty,
\end{equation}
uniformly with respect to $\alpha\ge0$, when $n$ is such that convergence at infinity  of the integral is guaranteed if $\alpha=0$, that is, if $n>\Re(a-b)$. 

\subsection{Expansion of \protectbold{{Q}_{-\frac12+i\tau}^\mu(x)} for large \protectbold{\tau}}\label{sec:QasyKum2}
For the representation of the conical $Q-$function in \eqref{eq:QKum01} we have $a=\frac12+\mu$, $b=\frac12-\mu$, $c=1$, $\omega=i\tau$. This gives
\begin{equation}\label{eq:QKum13}
\begin{array}{@{}r@{\;}c@{\;}l@{}}
e^{-\pi i\mu}{Q}_{-\frac12+i\tau}^\mu(x)&=&
\dsp{\sqrt{\pi/2}\, \alpha^{\mu+\frac12}\left(x^2-1\right)^{-\frac14}
\ \times}\\[8pt]
&&\dsp{\frac{1}{\Gamma\left(\frac12-\mu\right)}\int_0^\infty t^{-\mu-\frac12}(t+\alpha)^{-\mu-\frac12}e^{-\omega t}f(t)\,dt,}
\end{array}
\end{equation}
where $\alpha=\ln\frac{z+1}{z}$,  $z$ is defined in \eqref{eq:QKum02}, and
\begin{equation}\label{eq:QKum14}
f(t)=\left(\left(\frac{1-e^{-t}}{t}\right)\left(\frac{e^t-e^{-\alpha}}{t+\alpha}\right)\left(\frac{\alpha}{1-e^{-\alpha}}\right)\right)^{-\mu-\frac12}.
\end{equation}
By expanding $\dsp{f(t)=\sk f_k t^k}$ (observe that $f(0)=1$) we obtain
\begin{equation}\label{eq:QKum15}
\begin{array}{@{}r@{\;}c@{\;}l@{}}
e^{-\pi i\mu}{Q}_{-\frac12+i\tau}^\mu(x)&\sim&
\dsp{\sqrt{\pi/2}\, \alpha^{\mu+\frac12}
\left(x^2-1\right)^{-\frac14}\ \times}\\[8pt]
&&\dsp{\left(x+\sqrt{x^2-1}\right)^{-i\tau}\ \sn f_k \Phi_k,}
\end{array}
\end{equation}
where
\begin{equation}\label{eq:QKum16}
\Phi_k= \left(\tfrac12-\mu\right)_k \omega^{2\mu-k}U\left(\tfrac12+\mu,1+2\mu-k,\alpha \omega\right),
\quad \omega=i\tau.
\end{equation}

The first $f_k$ coefficients\footnote{More coefficients can be
obtained by using the Maple script available
at http://personales.unican.es/gila/fkcoeff.txt} are given by

\begin{equation}
\begin{array}{lcl}
f_0&=&1,\\[6pt]
f_1&=&\Frac{b}{2d}\left(2dz+d-2z\right),\\[6pt]
f_2&=&\Frac{b}{24d^2}\left(12z^2+12bz^2+d^2-12d^2z-12d^2z^2-24bdz^2\right.\\[6pt]
     &&\left.\,\,\,\,\,\,\,\,\,\, +12bd^2z+12bd^2z^2+3bd^2-12bdz\right),
\end{array}
\end{equation}
where $b=-\mu-\frac12$ and $d=z \alpha$.

The $\Phi_k$ can be written in terms of the Hankel functions. The first one is
\begin{equation}\label{eq:QKum17}
\Phi_0=-\tfrac12i\sqrt{\pi}(\tau/\alpha)^\mu e^{\frac12i\alpha\tau}H_\mu^{(2)}\left(\tfrac12\alpha\tau\right).
\end{equation}
For $\Phi_1$ we have
\begin{equation}\label{eq:QKum18}
\begin{array}{@{}r@{\;}c@{\;}l@{}}
\Phi_1=\dsp{-\frac{d}{d\omega}\Phi_0}
&=&
\dsp{\tfrac12\sqrt{\pi}(\tau/\alpha)^\mu e^{\frac12i\alpha\tau}} \ \times\\[8pt]
&&\left(
\left(\mu/\tau+\tfrac12i\alpha\right)H_\mu^{(2)}\left(\tfrac12\alpha\tau\right)
+\tfrac12\alpha {H_\mu^{(2)}}^\prime\left(\tfrac12\alpha\tau\right)\right)\\[8pt]
&=&
\dsp{\tfrac14\alpha\sqrt{\pi}(\tau/\alpha)^\mu e^{\frac12i\alpha\tau}} \left(iH_\mu^{(2)}\left(\tfrac12\alpha\tau\right)
+H_{\mu-1}^{(2)}\left(\tfrac12\alpha\tau\right)\right).
\end{array}
\end{equation}
The representation in terms of of the Hankel functions is convenient, because it is easy to separate real and imaginary parts by using
\begin{equation}\label{eq:QKum19}
H_\mu^{(2)}(z)=J_\mu(z)-iY_\mu(z).
\end{equation}

\section{Numerical tests of the expansions}\label{numxsmall}

In order to test the range of validity of the expansions, we have first compared the results
with the real part of the value obtained in the computation of ${Q}^0_{-\frac{1}{2}+i\tau}(x)$
using the command {\bf LegendreQ} in Maple with 50 digits.

\begin{figure}
\caption{Relative errors in the computation 
of $\widetilde{{Q}}^{0}_{-\frac{1}{2}+i\tau}(x)$  obtained by comparing the real part of \eqref{eq:QKumF06}   
with $0\le k \le 8$ against Maple. The results obtained for $\tau=0.1,\,5,\,20$ are shown.
\label{Fig2}}
\begin{center}
\epsfxsize=12.5cm \epsfbox{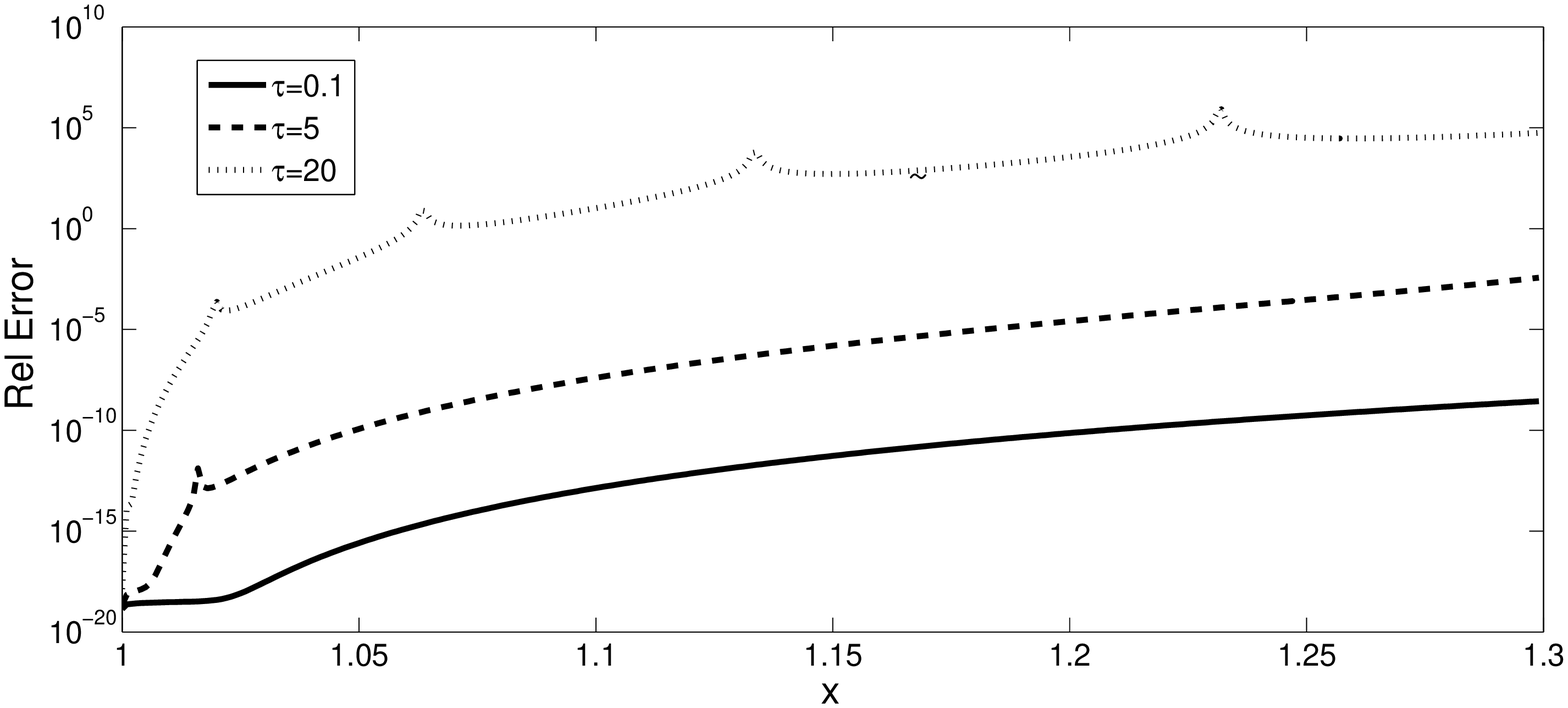}
\end{center}
\end{figure}

\begin{figure}
\caption{Relative errors in the computation 
of $\widetilde{{Q}}^{0}_{-\frac{1}{2}+i\tau}(x)$  obtained by comparing the real part of \eqref{eq:hyp1}   
with $0\le k \le 8$ against Maple. The results obtained for $\tau=0.1,\,5,\,20$ are shown.
\label{Fig3}}
\begin{center}
\epsfxsize=12.5cm \epsfbox{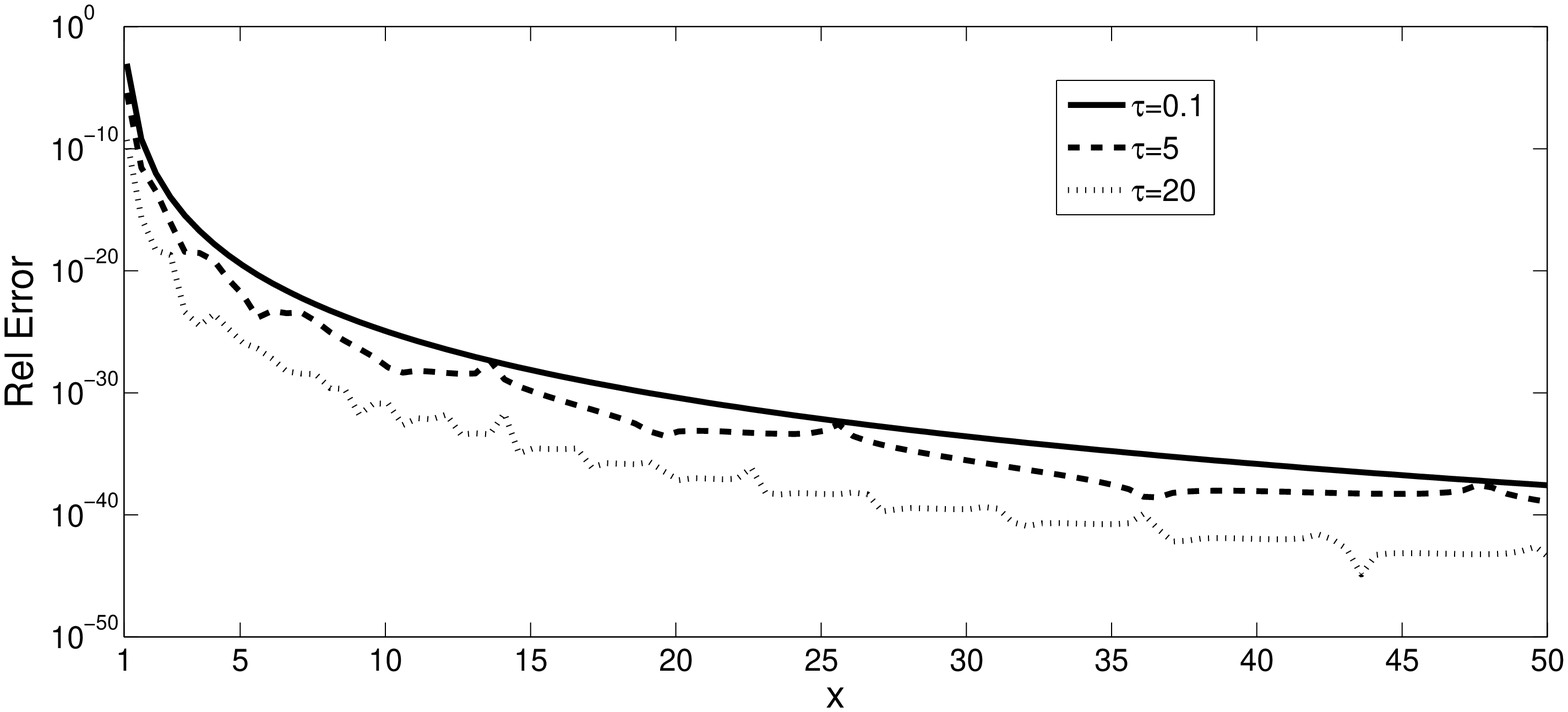}
\end{center}
\end{figure}

Fig.~\ref{Fig2} shows, as a function of $x$ and for three different values of $\tau$, the accuracy obtained
in the computation of  $\widetilde{{Q}}^{0}_{-\frac{1}{2}+i\tau}(x)$ using the power series
expansion given in \eqref{eq:QKumF06}. We have used the expansion with $0 \le k \le 8$. As can be seen,
the expansion provides full accuracy for values of $x$ very close to 1 but the accuracy worsens as $\tau$ increases. 
As an example, the number of terms of the expansion has to be increased to more than $100$ in order
to obtain an accuracy better $10^{-14}$ when computing  $\widetilde{{Q}}^{0}_{-\frac{1}{2}+i200}(1.05)$. 
The same test has been performed using the power series of \eqref{eq:hyp1} for the
computation of   $\widetilde{{Q}}^{0}_{-\frac{1}{2}+i\tau}(x)$. The results are shown in Fig.~\ref{Fig3}.
In this case, as expected the expansion behaves better as $\tau$ increases but the accuracy worsen as $x$ approaches
to $1$.

\begin{figure}
\caption{Relative errors in the computation 
of $\widetilde{{Q}}^{0}_{-\frac{1}{2}+i\tau}(x)$  obtained by comparing the real part of \eqref{eq:QKum15}   
with $n=8$ against Maple. The results obtained for $x=1.1,\,100$ are shown.
\label{Fig4}}
\begin{center}
\epsfxsize=12.5cm \epsfbox{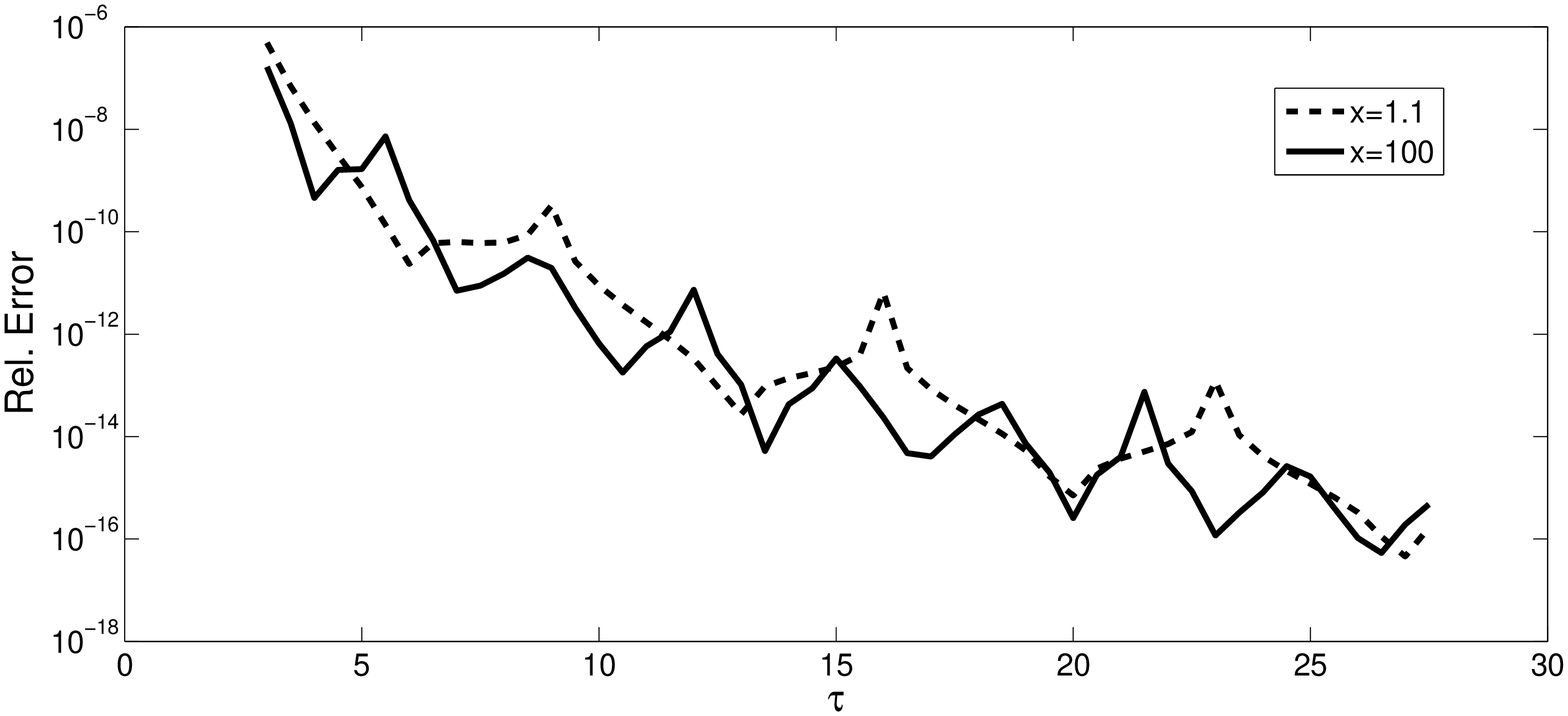}
\end{center}
\end{figure}

A test for the expansion in terms of the Kummer $U$-function\footnote{A Maple implementation is available
at http://personales.unican.es/gila/QKummer.txt  } is shown in
 Fig.~\ref{Fig4}: the relative errors obtained in the computation of 
$\widetilde{{Q}}^{0}_{-\frac{1}{2}+i\tau}(x)$  using the expansion \eqref{eq:QKum15}
 with $0\le k \le 8$ in the computation are plotted. 
 The  functions $\Phi_k$ in \eqref{eq:QKum15} are computed by using the recurrence relation \eqref{eq:QKum10}
with starting values those given in \eqref{eq:QKum17} and \eqref{eq:QKum18}.
As can be seen, the expansion provides an accuracy better than $10^{-8}$ (single
precision) even for values of the parameter $\tau$ lower than $5$. Two values of the argument $x$ are chosen ($x=1.1,\,100$) 
in order to illustrate the validity of the expansion for large and small (larger than $1$) values of $x$. 

  A test independent of the Maple procedure {\bf LegendreQ} has also been considered. It 
consists in testing that the values computed using the expansion 
are consistent with the three-term recurrence relation \eqref{TTRR2} for $\widetilde{{Q}}^{m}_{-\frac{1}{2}+i\tau}(x)$. 
We first compute  $\widetilde{{Q}}^{0}_{-\frac{1}{2}+i\tau}(x)$, $\widetilde{{Q}}^{1}_{-\frac{1}{2}+i\tau}(x)$
and  $\widetilde{{Q}}^{2}_{-\frac{1}{2}+i\tau}(x)$ using the expansion and check
 \eqref{TTRR2} written in the form

\begin{equation}
\Frac{\left(2mx/\sqrt{x^2-1}\right)\widetilde{{Q}}^{m}_{-\frac{1}{2}+i\tau}(x)-\left((m-\tfrac12)^2+\tau^2\right)\widetilde{Q}_{-\frac12+i\tau}^{m-1}(x)}
{\widetilde{{Q}}^{m+1}_{-\frac{1}{2}+i\tau}(x) }=1.      
\label{errRR}
\end{equation}

\begin{figure}
\caption{Test of the recurrence relation given in \eqref{errRR} for $m=1$ and 
$x=1.1,\,100$ using the asymptotic expansion  \eqref{eq:QKum15}   
with $n=8$. 
\label{Fig5}}
\begin{center}
\epsfxsize=12.5cm \epsfbox{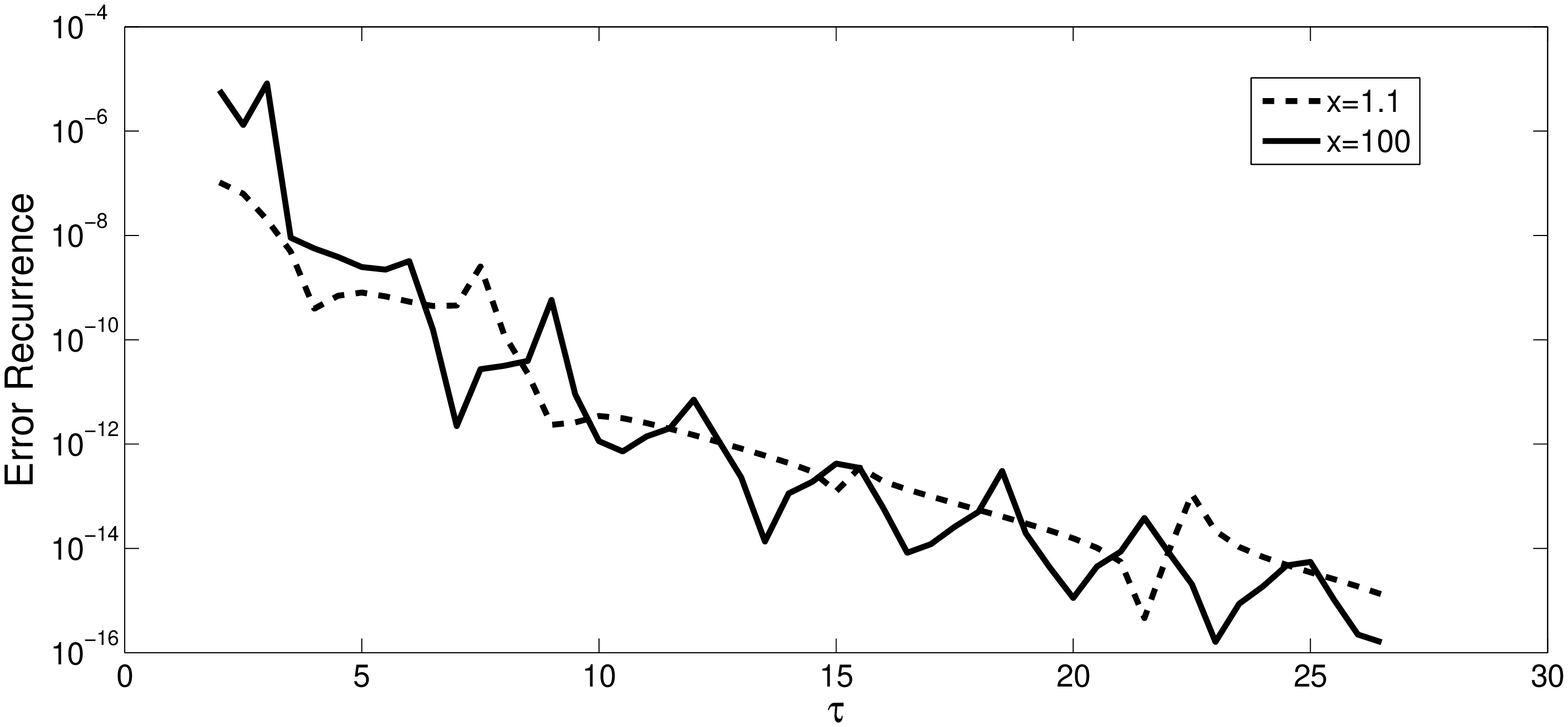}
\end{center}
\end{figure}

 The deviations from $1$ of the left-hand side of \eqref{errRR} for $m=1$ are shown in Fig.~\ref{Fig5}. The results are consistent
with those obtained when testing the computed value of  $\widetilde{{Q}}^{0}_{-\frac{1}{2}+i\tau}(x)$.  

\begin{figure}
\caption{Test of the recurrence relation given in \eqref{errRR} for $\tau$ fixed ($\tau=50$)
and  $x=1.00001,\,10,\,100,\,500$ using the asymptotic expansion  \eqref{eq:QKum15}   
with $n=8$. 
\label{Fig6}}
\begin{center}
\epsfxsize=12.5cm \epsfbox{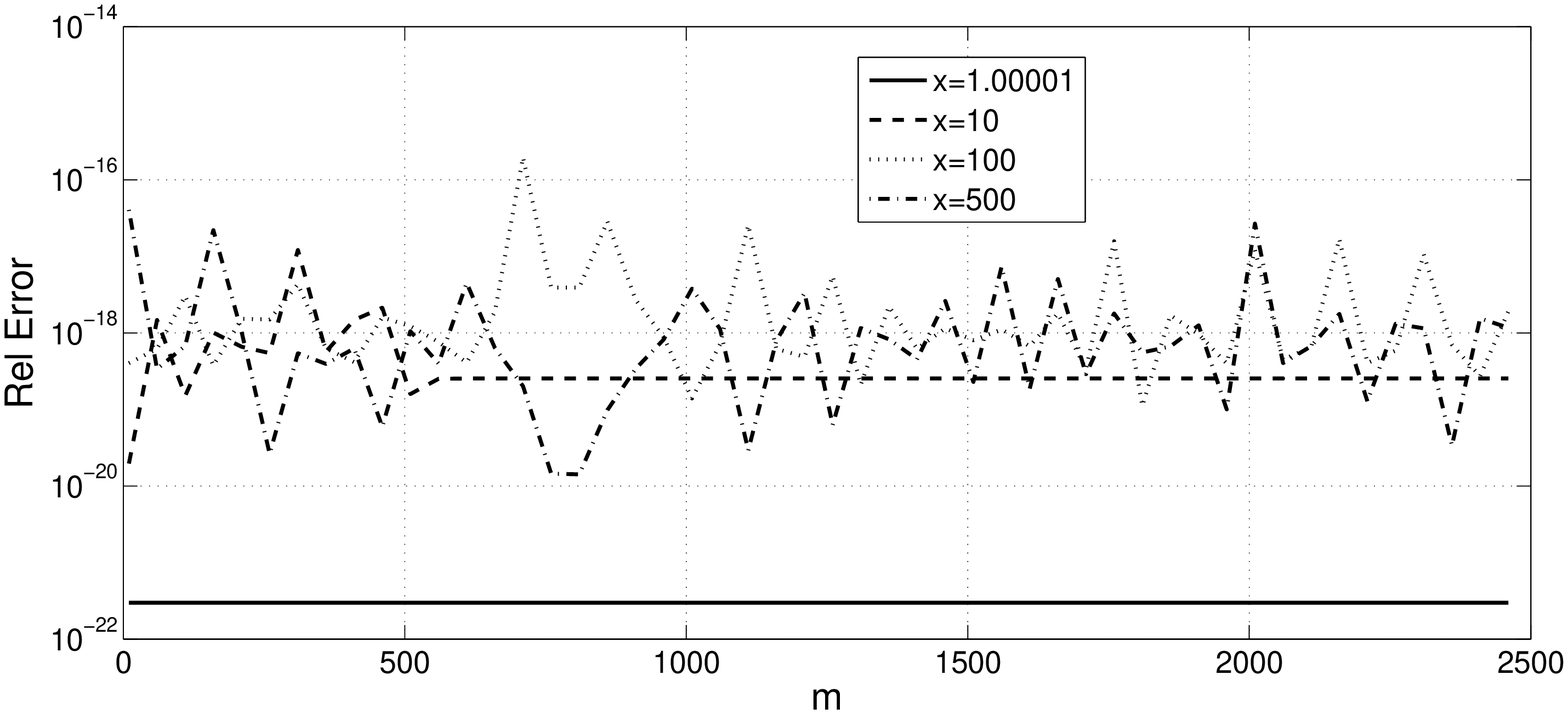}
\end{center}
\end{figure}

Another test involving the application of the three-term recurrence relation can be seen in Fig.~\ref{Fig6}. In this figure
we show as a function of $m$, the relative error obtained in the computation of $\widetilde{{Q}}^{m}_{-\frac{1}{2}+i\tau}(x)$ by 
applying the recurrence relation starting from $\widetilde{{Q}}^{0}_{-\frac{1}{2}+i\tau}(x)$ and $\widetilde{{Q}}^{1}_{-\frac{1}{2}+i\tau}(x)$.
The value of $\tau$ has been fixed to $50$ and four different values of $x$ ($x=1.00001,\,10,\,100,\,500$) have been considered
in the computations. As can be seen, the error is quite uniform even when very large values of $m$ are considered. This illustrates
the numerical stability of the forward application of the recursion for the $\widetilde{{Q}}^{m}_{-\frac{1}{2}+i\tau}(x)$, as expected from the
stability analysis of Sect.~\ref{RR}. Other (large) values of the parameter $\tau$ have been also tested, providing similar results.

\section{Computational scheme}

From the results obtained in the previous section, a stable computational scheme for evaluating the function
$\widetilde{{Q}}^{m}_{-\frac{1}{2}+i\tau}(x) \equiv 
\Re\left\{e^{-i\pi m} {{Q}}^{m}_{-\frac{1}{2}+i\tau}(x)\right\}$, a numerically satisfactory companion
of the function ${\rm P}^{m}_{-\frac{1}{2}+i\tau}(x)$ when $x>1$, emerges:
 
\begin{description}
\item{When $m=0,\,1$:}
\begin{enumerate}
\item{For} $x$ close to $1$ and small/moderate values of $\tau$, compute $\widetilde{{Q}}^{m}_{-\frac{1}{2}+i\tau}(x)$
 using the power series for the hypergeometric representation given in Sect.~\ref{sec:pow2}.
\item{For} $x$ close to $1$ and moderate/large values of $\tau$, compute $\widetilde{{Q}}^{m}_{-\frac{1}{2}+i\tau}(x)$
 using the expansion for large $\tau$ given in Sect.~\ref{sec:QasyKum2}.
\item{For} other values of $x$, compute $\widetilde{{Q}}^{m}_{-\frac{1}{2}+i\tau}(x)$ using the power series 
for the hypergeometric representation given in Sect.~\ref{sec:pow3}.
\end{enumerate}

\item{When $m \ge 2$:}\\ compute $\widetilde{{Q}}^{0}_{-\frac{1}{2}+i\tau}(x)$ and $\widetilde{{Q}}^{1}_{-\frac{1}{2}+i\tau}(x)$
using the previous scheme and then use the recursion 

\begin{equation}
\label{TTRRQ}
\widetilde{Q}_{-\frac12+i\tau}^{m +1}(x)-\Frac{2m x}{\sqrt{x^2-1}}
\widetilde{Q}_{-\frac12+i\tau}^{m}(x) +
\left((m-\tfrac12)^2+\tau^2\right)\widetilde{Q}_{-\frac12+i\tau}^{m-1}(x)  =0
\end{equation}
in the direction of increasing $m$.   
\end{description}

Numerical tests implemented in Maple show that an accuracy better than $10^{-14}$ can be obtained
for computing the initial values of the recurrence relation when 
\begin{description}
\item[a.] 
the hypergeometric representation 
of Sect.~\ref{sec:pow2} is used when $x<1.1$ and $\tau<10$;
\item[b.]
the expansion for large $\tau$ of Sect.~\ref{sec:QasyKum2} is used  when $x<1.1$ and $\tau\ge 10$; 
\item[c.]
 the hypergeometric representation
of Sect.~\ref{sec:pow3} is used in other cases.
\end{description}

For higher accuracies, the regions of application of the methods may need some adjustments. 

\section*{Acknowledgements}
The authors thank the reviewers for useful comments.
This work was supported by  {\emph{Ministerio de Ciencia e Innovaci\'on}}, 
project MTM2009-11686 and {\emph{Ministerio de Econom\'{\i}a y Competitividad}},
project MTM2012-34787. NMT thanks CWI, Amsterdam, for scientific support.

\bibliographystyle{plain}    
\bibliography{Conical} 

\begin{thebibliography}{1}

\bibitem{Bateman:1953:HTF}
H.~Bateman.
\newblock {\em Higher Transcendental Functions}, volume~I.
\newblock McGraw-Hill, 1953.

\bibitem{dunster:1991:cfo}
T.~M. Dunster.
\newblock Conical functions with one or both parameters large.
\newblock {\em Proc. Roy. Soc. Edinburgh Sect. A}, 119(3-4):311--327, 1991.

\bibitem{Duns:2010:Con}
T.~M. Dunster.
\newblock Legendre and related functions.
\newblock In {\em N{IST} handbook of mathematical functions}, pages 351--381.
  Cambridge University Press, New York, 2010.

\bibitem{Gil:2007:NSF}
A.~Gil, J.~Segura, and N.~M. Temme.
\newblock {\em Numerical methods for special functions}.
\newblock Society for Industrial and Applied Mathematics (SIAM), Philadelphia,
  PA, 2007.

\bibitem{gil:2012:cpc}
A.~Gil, J.~Segura, and N.~M. Temme.
\newblock An improved algorithm and a fortran 90 module for computing the
  conical function ${P}^m_{-1/2+i\tau}(x)$.
\newblock {\em Comput Phys Commun}, 183:794--799, 2012.

\bibitem{gil:2009:con}
A.~Gil, J.~Segura, and N.M. Temme.
\newblock Computing the conical function {$P^\mu_{-1/2+i\tau}(x)$}.
\newblock {\em SIAM J. Sci. Comput.}, 31(3):1716--1741, 2009.

\bibitem{Temme:1983:NCC}
N.~M. Temme.
\newblock The numerical computation of the confluent hypergeometric function
  ${U}(a,b,z)$.
\newblock {\em Numer. Math.}, 41:63--82, 1983.

\bibitem{Temme:1996:SFA}
N.~M. Temme.
\newblock {\em Special functions: {A}n introduction to the classical functions
  of mathematical physics}.
\newblock A Wiley-Interscience Publication. John Wiley \& Sons Inc., New York,
  1996.

\bibitem{thebault:2006:geo}
E.~Thebault.
\newblock Revised cap harmonic analysis (r-scha): Validation and properties.
\newblock {\em J. Geophys. Research}, 111, 2006.

\end{thebibliography}

\end{document}